\documentclass[12pt,reqno]{amsart} 
\usepackage{amsmath,amssymb,amsthm,mathrsfs}
\usepackage[english]{babel}
\usepackage[utf8]{inputenc}
\usepackage{xcolor}

\newtheorem{theorem}{Theorem}[section]
\newtheorem{proposition}[theorem]{Proposition}

\theoremstyle{definition}
\newtheorem{definition}[theorem]{Definition}
\newtheorem{example}[theorem]{Example}
\newtheorem{remark}[theorem]{Remark}

\usepackage{multirow} 
\usepackage{color}

\usepackage{graphicx}
\usepackage[left=3cm,right=3cm]{geometry}
\usepackage[T1]{fontenc}

\usepackage{cancel,soul,ulem}

\def\r{\mathbb R}
\def\h{\mathbb H}
\def\s{\mathbb S}

\begin{document}

\title{Conformal trajectories in $3$-dimensional space forms}
\author{Rafael L\'opez, Marian Ioan Munteanu}
\address{ Departamento de Geometr\'{\i}a y Topolog\'{\i}a\\  Universidad de Granada. 18071 Granada, Spain}
\email{rcamino@ugr.es}
\address{University 'Al. I. Cuza' of Iasi, Faculty of Mathematics, Bd. Carol I, no. 11, 700506 Iasi, Romania}
\email{marian.ioan.munteanu@gmail.com}

\begin{abstract}
We introduce the notion of conformal trajectories in three-dimensional Riemannian manifolds $M^3$. Given a conformal vector field   $V\in\mathfrak{X}(M^3)$, a conformal trajectory  of $V$ is a regular curve $\gamma$ in $M^3$  satisfying  $\nabla_{\gamma'}\gamma'=q\, V\times\gamma'$, for some  fixed non-zero constant $q\in\r$.   In this paper, we study conformal trajectories in the space forms  $\r^3$, $\s^3$ and $\h^3$. For (non-Killing) conformal  vector fields in $\s^3$ 
(respectively in $\h^3$), we prove that conformal trajectories have constant curvature and its   torsion is a linear combination of trigonometric (respectively hyperbolic) functions on the arc-length parameter.  In the case of Euclidean space $\r^3$, we obtain the same result for the radial vector field and characterizing  all conformal trajectories.
 \end{abstract}


\maketitle

\section{Introduction and motivation}

 In a $n$-dimensional Riemannian manifold $(M^n,\langle~,~\rangle)$, a magnetic field is a closed $2$-form $F$. A curve $\gamma\colon I\to M^n$ is called a magnetic trajectory of $F$ if $\gamma$ satisfies the Lorentz equation $\nabla_{\gamma'}\gamma'=\Phi(\gamma')$, where $\Phi$ is the $(1,1)$-tensor field on $M$ defined by $\langle\Phi X,Y\rangle=F(X,Y)$, $X,Y\in \mathfrak{X}(M^n)$.  The notion of magnetic trajectory extends that of geodesic when $F=0$. 

Suppose that the dimension of $M^n$ is $n=3$.   In such a case, the study of  magnetic trajectories can be introduced with  a different approach. This is due  because $2$-forms and vector fields are identified via the Hodge star operator $*$ and the volume form $\omega$ of $M^3$. In particular, magnetic fields are vector fields with zero divergence and the Lorentz force can be expressed in terms of the cross product $\times$ of $M^3$.   Killing vector fields are examples of magnetic fields.  Each Killing vector field   defines a magnetic field $F_V$ by means of $F_V=i_V\omega$. Then the   Lorentz force $\Phi$ of $F_V$ is given by  
$$\Phi(X)=V\times X.$$
Thus a magnetic trajectory $\gamma$ of $V$ is characterized by the identity
\[
\nabla_{\gamma'}\gamma'=V\times\gamma'.
\]
Looking this identity, the property that $V$ is magnetic can be ignored and it only uses the cross product of $M^3$. This observation allows us to generalize  the notion of   magnetic curve by replacing a Killing vector field by a conformal  vector field.   A vector field   $V\in\mathfrak{X}(M^n)$ is said to be {\it conformal} if 
$$\mathcal{L}_V\langle~,~\rangle =\lambda \langle~,~\rangle,$$
 where $\mathcal{L}$ is the Lie derivative and $\lambda\in C^\infty(M^n)$ is a function on $M^n$. Equivalently, we have  
\begin{equation}\label{conformal}
\langle \nabla_XV,Y\rangle+\langle X,\nabla_YV\rangle=\lambda \langle X,Y\rangle,\quad X,Y\in\mathfrak{X}(M^n).
\end{equation}
In the literature, it is also used the terminology of conformal Killing vector field. If $\lambda=0$, then $V$ is just a Killing vector field. Those ones that are not Killing   will be caller proper conformal vector fields. 

\begin{definition}   Let $V\in\mathfrak{X}(M^3)$ be a conformal  vector field. A {\it conformal trajectory} of $V$ is a regular curve $\gamma\colon I\subset\r\to M^3$ satisfying the equation
\begin{equation}\label{eq1}
\nabla_{\gamma'}\gamma'=q\,V\times\gamma',
\end{equation}
where $q\in\r$, $q\not=0$,  is fixed.
\end{definition}

The number $q$ in \eqref{eq1} is due because conformal trajectories cannot be rescaled, such as it occurs in  the case of  geodesics. Indeed, if $\gamma$ is a conformal trajectory, then 
$$\frac{d}{ds} \langle \gamma'(s),\gamma'(s)\rangle=2\langle \nabla_{\gamma'(s)}\gamma'(s),\gamma'(s)\rangle=2 \langle V\times\gamma',\gamma'\rangle=0.$$
As a consequence the speed $|\gamma'(s)|$ of $\gamma$ is constant. 

 The purpose of this paper is to study  conformal trajectories in the space forms $\r^3$, $\s^3$ and $\h^3$.  For magnetic  fields, the study of magnetic trajectories in space forms have been developed in \cite{dm1,dm2,ei,ke,mu,mn,wa}. 
 In each of these space forms, we will consider particular examples of proper conformal  vector fields. To be precise, 
 \begin{enumerate}
 \item In Euclidean space $\r^3$, let $V$ be the radial vector field $V(x,y,z)=x\partial_x+y\partial_y+z\partial_z$.
 \item In the unit sphere $\s^3$,   proper conformal vector fields are generated by $V(p)=\vec{a}-\langle \vec{a},p\rangle p$, where $\vec{a}\in\r^4$.  
  \item In the hyperbolic space $\h^3$, consider  the Lorentzian model  where $\h^3$ is viewed as a quadric of $\r^4$. The proper conformal vector fields are generated by  $V(p)=\vec{a}+\langle\vec{a},p\rangle p$, where $\vec{a}\in\r^4_1$.

 \end{enumerate}
 In all cases, when we express Eq. \eqref{eq1} in coordinates, the equation turns on a second order differential equation whose explicit solutions are difficult to find in general. However, we can summarize a property of the conformal trajectories. 
 
 \begin{theorem} Let $M^3$ be a space form and let $V\in\mathfrak{X}(M^3)$ be the conformal vector field described above. Let $\gamma\colon I\subset\r\to M^3$ be a regular (non-geodesic) curve parametrized by arc-length $\gamma=\gamma(s)$. If $\gamma$ is a conformal trajectory of $V$, then its curvature $\kappa$ is a non-zero constant and its torsion $\tau$ is:
 \begin{enumerate}
 \item in  $\r^3$,  $\tau(s)=a_1 s+a_2$;
  \item in  $\s^3$, $\tau(s)=a_1 \sin s+a_2\cos s $;
  \item in   $\h^3$, $\tau(s)=a_1 \sinh s+a_2\cosh s$,
 \end{enumerate}
 where $a_1,a_2\in\r$.
 \end{theorem}
 
 The proof of this theorem will be done in Sects. \ref{s3}, \ref{s4} and \ref{s5} in each one of the space forms $\r^3$, $\s^3$ and $\h^3$, respectively.  In Sect. \ref{s2} we present all conformal vector fields in space forms and we calculate the curvature and the torsion of conformal trajectories for  proper conformal vector fields.
 
\section{Conformal vector fields in space forms}
\label{s2}

In $3$-dimensional space forms, the space of conformal vector fields is $10$-dimensional, 
from which $6$ correspond to Killing vector fields. We describe in each space form these vector fields. General references are \cite{di,st}.
\begin{enumerate}
\item Euclidean space $\r^3$. Consider $(x,y,z)$ canonical coordinates in $\r^3$ and let 
$\langle~,~\rangle=dx^2+dy^2+dz^2$ be the Euclidean metric. The Killing vector fields are generated by the translations $\partial_x$, $\partial_y$ and $\partial_z$ and the rotations $-y\partial_x+x\partial_y$, $-z\partial_x+x\partial_z$ and $-z\partial_y+y\partial_z$. The proper conformal vector fields are generated by the radial vector field
\begin{equation}\label{v1}
V=x\partial_x+y\partial_y+z\partial_z,
\end{equation}
and the vector fields
\begin{equation}\label{w1}
W(p)=\frac{|p|^2}{2}\vec{a}-\langle \vec{a},p\rangle p,\quad p\in\r^3,
\end{equation}
where $\vec{a}\in\r^3$.  The vector fields $V$ and $W$ are conformal because 
$\mathcal{L}_V \langle~,~\rangle=2\langle~,~\rangle$ and  
$\mathcal{L}_{W}\langle~,~\rangle=-2\langle\vec{a},p\rangle \langle~,~\rangle$.

\item Sphere $\s^3$. Let $\langle~,~\rangle$ be the Euclidean metric on $\s^3$. 
Let $(x,y,z,t)$ be canonical coordinates in $\r^4$ and let identify a   vector field with its coordinates $(x,y,z,t)$ with respect to the canonical basis  $\{\partial_x,\partial_y,\partial_z,\partial_t\}$ of $\r^4$. Then the Killing vector fields are generated by 
\begin{eqnarray*}
&K_1=(-y,x,-t,z),\quad &K_2=(-y,x,t,-z),\\
&K_3=(-z,-t,x,y,)\quad &K_4=(-z,t,x,-y),\\
&K_5=(-t,z,-y,x),\quad &K_6=(-t,-z,y,x).
\end{eqnarray*}
Proper conformal vector fields are the projections on $\s^3$ of the analogous  vector fields $W$ given in \eqref{w1} but defined in $\r^4$. Thus, there is a description of these vector fields free of coordinates, namely,  
\begin{equation}\label{v2}
V(p)=\vec{a}-\langle \vec{a},p\rangle p,\quad p\in\s^3,
\end{equation}
 where $\vec{a}\in\r^4$. Here $\mathcal{L}_{V}\langle~,~\rangle=-2\langle\vec{a},p\rangle 
 \langle~,~\rangle$.
 \item Hyperbolic space $\h^3$. Consider the Lorentzian model, that is, the Euclidean space  
 $\r^4$ with canonical coordinates $(x,y,z,t)$, endowed with  the Lorentzian metric   
 $\langle~,~\rangle=dx^2+dy^2+dz^2-dt^2$. 
 Then $\h^3$ is the quadric $\h^3=\{p\in\r^4:\langle p,p\rangle=-1,~t>0\}$.  As in $\s^3$, we identify a vector field with its coordinates with respect to the canonical basis of tangent vectors. The   Killing vector fields are generated by 
 \begin{eqnarray*}
&K_1=(-y,x,t,z)\quad &K_2=(y,-x,t.z)\\
&K_3=(-z,t,x,y)\quad &K_4=(z,t,-x,y)\\
&K_5=(t,-z,y,x)\quad &K_6=(t,z,-y,x).
\end{eqnarray*}
 Similarly as in the case of the sphere $\s^3$, the space of   proper conformal vector fields are generated by 
\begin{equation}\label{v3}
V(p)=\vec{a}+\langle \vec{a},p\rangle p,\quad p\in\h^3,
\end{equation}
 where $\vec{a}\in\r^4$. In this case $\mathcal{L}_{W}\langle~,~\rangle=2\langle\vec{a},p\rangle \langle~,~\rangle$.
\end{enumerate}

\medskip

Let $V\in\mathfrak{X}(M^3)$ be a conformal  vector field with $\mathcal{L}_V\langle~,~\rangle=\lambda \langle~,~\rangle$. Let $\gamma\colon I\to M^3$ be a regular curve parametrized by arc-length. If $\gamma$ is a geodesic, then it  is immediate the following result. 

\begin{proposition}\label{pr1}
 A geodesic  $\gamma\colon I\to M^3$ is a conformal trajectory of $V$ if and only if $V(\gamma(s))$ and $\gamma'(s)$ are collinear for all $s\in I$.
\end{proposition}

 From now on, suppose that $\gamma$ is not a geodesic. Then we can assign to $\gamma$ its Frenet frame $\{T,N,B\}$, where $T(s)=\gamma'(s)$ and $T\times N=B$. 
  The Frenet equations are given by
\begin{equation*}
\begin{split}
\nabla_TT&=\kappa N,\\
\nabla_TN&=-\kappa T+\tau B,\\
\nabla_TB&=-\tau B,
\end{split}
\end{equation*}
where $\kappa$ and $\tau$ are the curvature and the torsion  of $\gamma$, respectively. Let denote $V(s)=V(\gamma(s))$ 

In the following result, we study the tangent part of $V$ along $\gamma$ for those   conformal vector fields $V$ satisfying that $\nabla_TV$ is proportional to $T$ along $\gamma$.

\begin{proposition}\label{pr2}
   Let $\gamma\colon I\to M^3$ be a   curve parametrized by arc-length. If $\gamma$ is a conformal trajectory of $V$ then
\begin{equation}\label{aa}
 \langle  V(s),\gamma'(s)\rangle=\frac12\int^s_{}\lambda(\gamma(t))dt.
\end{equation}
If, in addition, $\nabla_TV$ is proportional to $T$ for all $s\in I$, then the curvature $\kappa$ is constant, $\kappa(s)=\kappa_0$ and the torsion is 
\begin{equation}\label{tt}
\tau(s)=q\langle V(s),\gamma'(s)\rangle.
\end{equation}
\end{proposition}

\begin{proof}
We know   $\nabla_{\gamma'}\gamma'=\nabla_TT=\kappa N$. We now compute the right hand-side of \eqref{eq1}. For this, let us write
\begin{equation}\label{2}
V(s)=a(s) T(s)+b(s) N(s)+c(s) B(s),
\end{equation}
where $a,b,c\in C^\infty(I)$. Then $V\times\gamma'=V\times T=cN-bB$. Thus Eq. \eqref{eq1} becomes 
$$\kappa N=q( cN-bB).$$
This implies $b=0$ and $\kappa=qc$. In particular, we have $c\not=0$. 
Using this, we obtain
\begin{equation}\label{vv}
\nabla_{\gamma'}V=a'T+\kappa(a-\frac{\tau}{q})N+\frac{\kappa'}{q}B.
\end{equation}
Multiplying by $T$ and using \eqref{conformal} because $V$ is a conformal  vector field, we have
\[
\lambda=\lambda \langle \gamma',\gamma'\rangle= 2\langle \nabla_{\gamma'}V,\gamma'\rangle=2a'.
\]
This proves \eqref{aa}. If now we assume that $\nabla_TV$ is proportional to $T$, 
then \eqref{vv} yields $\kappa'(s)=0$ for all $s\in I$ and $\tau(s)=qa(s)$. 
Thus $\kappa$ is a nonzero constant and the torsion is given by \eqref{tt} 
because of \eqref{aa}.
\end{proof}

We can use this proposition for the proper conformal vector field in space forms described previously. 

 \begin{theorem}\label{t0}
  Let $\gamma\colon I\to M^3$ be a   curve parametrized by arc-length, where $M^3$ is a space form $\r^3$, $\s^3$ or $\h^3$. Suppose that $\gamma$ is not a geodesic. 
 \begin{enumerate}
 \item Case $\r^3$. Let $V$ be the radial vector field in $\r^3$ defined in \eqref{v1}. If $\gamma$ is a conformal trajectory of $V$, then the curvature of $\gamma$ is constant and its torsion is 
 $$\tau(s)=q\langle \gamma(s),\gamma'(s)\rangle.$$
 \item Case $\s^3$. Let $V$ be a proper conformal vector field defined in \eqref{v2}. If $\gamma$ is a conformal trajectory of $V$, then the curvature of $\gamma$ is constant and its torsion is 
 $$\tau(s)=q\langle\vec{a},\gamma'(s)\rangle.$$
 \item Case $\h^3$.  Let $V$ be a proper conformal vector field defined in \eqref{v3}. If $\gamma$ is a conformal trajectory of $V$, then the curvature of $\gamma$ is constant and its torsion is 
 $$\tau(s)=q\langle\vec{a},\gamma'(s)\rangle.$$
 \end{enumerate}
 \end{theorem}
 
 \begin{proof}In all cases, we have $\nabla_TV$ is proportional to $T$, hence we can apply Prop. \ref{pr2}.
 \end{proof}

 \section{The Euclidean space $\r^3$: the radial vector field}\label{s3}
 In this section we study the conformal trajectories of the radial vector field $V$ defined in \eqref{v1}.   Notice that the conformal transformations of $V$ are dilations from the origin. By Prop. \ref{pr1},  a straight-line   $\gamma(s)=p+s\vec{v}$  is a conformal trajectory of $V$ if and only if $p$ and $\vec{v}$ are collinear, that   $\gamma$ is the straight-line acrosses the origin of $\r^3$. 

From now, we will assume that $\gamma$ is not a straight-line. Suppose that $\gamma$ is parametrized by arc-length. The left hand-side of \eqref{eq1} is $\nabla_{\gamma'}\gamma'=\gamma''$. Thus $\gamma$ is a conformal trajectory if and only if
\begin{equation}\label{eq2}
\gamma''=q\, \gamma\times\gamma'.
\end{equation}

\begin{theorem} \label{t1}
If $\gamma$ is a  conformal trajectory   of the radial vector field  \eqref{v1}, 
then $\gamma$  is a curve with constant curvature $\kappa(s)=\kappa_0$, and its 
torsion is an affine function in the parameter of the curve. Moreover, 
\begin{equation}\label{3}
\gamma(s)=(s+a_0)\gamma'(s)+\frac{1}{q}\gamma'(s)\times\gamma''(s).
\end{equation}
\end{theorem}

\begin{proof}

 Notice that $\nabla_TV=\nabla_{\gamma'}V=\gamma'=T$. Thus we are under the assumptions of Prop. \ref{pr2}, in particular,   the curvature $\kappa$  is constant and \eqref{aa} gives 
 $$\langle \gamma(s),\gamma'(s)\rangle=\int^s \ ds= s+a_0,$$
where $a_0\in\r$ is an integration constant. Moreover, from \eqref{tt}, we obtain $\tau(s)=q(s+a_0)$.  Using the same notation as in the proof of Prop.~\ref{pr2}, equation \eqref{vv}  implies     $c(s)=c_0\not=0$ is constant and 
$\kappa(s)=qc_0$.   Notice that $cB=\frac{\kappa}{q}T\times N=\frac{1}{q}\gamma'\times\gamma''$. This yields \eqref{3}.
 \end{proof}
 
Equation \eqref{3} implies that the radial vector of $\gamma$ lies in the rectifying plane at each point. This class of curves were studied by B-Y. Chen in \cite{ch}. However  identity \eqref{3} does not characterize conformal trajectories of the radial vector field. This is because   the curvature and torsion determine a spatial curve of $\r^3$, up to rigid motions of $\r^3$. However,  identity \eqref{eq1} depends on the position vector of the curve. 

A consequence of  \eqref{3} is
 $$ |\gamma(s)|^2= (s+a_0)^2+c_0^2.$$
This allows to give the converse of Thm. \ref{t1}, obtaining a characterization of conformal trajectories of the radial vector field. 

\begin{theorem} Let $\gamma\colon I\subset\r\to\r^3$ be a curve parametrized by arc-length. Suppose that the curvature of $\gamma$ is constant, $\kappa(s)=\kappa_0$ and that there are $a_0,a_1\in\r$, $a_1>0$ such that 
\begin{equation}\label{4}
 |\gamma(s)|^2= (s+a_0)^2+a_1^2.
 \end{equation}
 Then $\gamma$ is a conformal trajectory of $V$ for some $q\in\r$.
 \end{theorem}
 \begin{proof}
 Differentiating twice \eqref{4}, we have $\langle  \gamma,T\rangle=s+a_0$ and $\langle \gamma,N\rangle=0$. By expressing $\gamma$ in terms of the Frenet frame, and using \eqref{4}, we obtain 
 $$\gamma(s)=(s+a_0)T+a_1 B.$$
 Then $\gamma\times\gamma'=a_1 N$. On the other hand, the first Frenet equation gives $ \gamma''=k_0 N$. Then taking $q=\frac{\kappa_0}{a_1}$, we have identity \eqref{eq2}. 
 \end{proof}

\begin{remark} 
For a curve $\gamma$ in Euclidean space $\r^3$ parametrized by arc-length $\gamma=\gamma(s)$, 
the identity \eqref{4} is  a characterization to be a  rectifying curve: see \cite[Theorem~1]{ch}.
\end{remark}

\begin{figure}[hbtp] 
\begin{center}
\includegraphics[width=.4\textwidth]{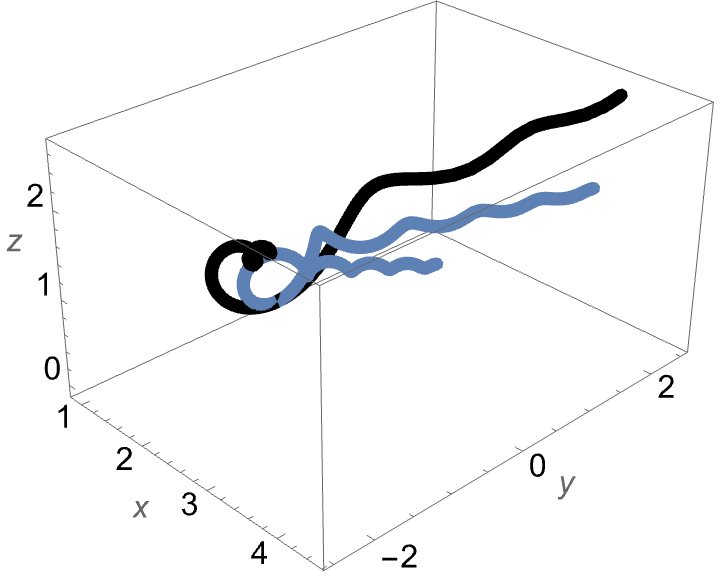} \qquad \includegraphics[width=.4\textwidth]{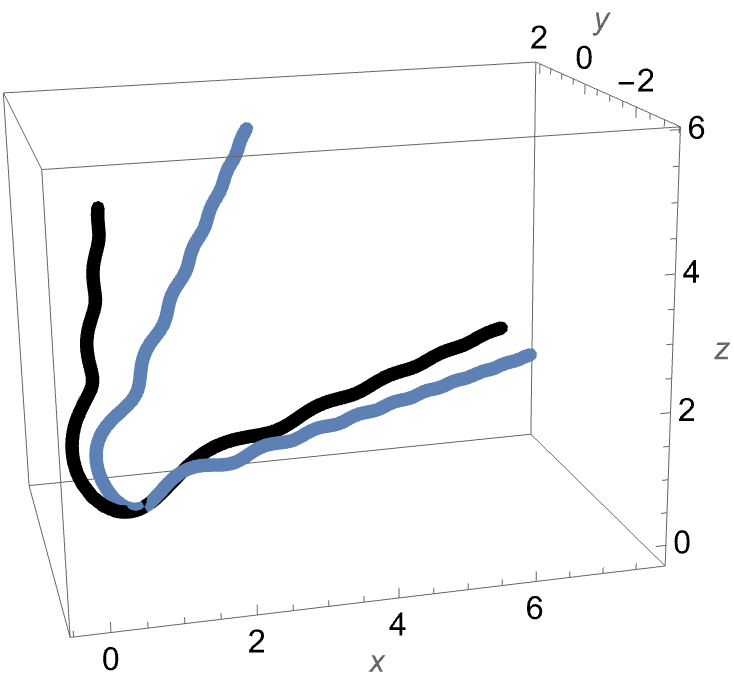}
\end{center}
\caption{Conformal trajectories of the radial vector field \eqref{v1} with different values for $q$: $q=1$ (black) and $q=2$ (grey). Left: initial conditions \eqref{ini1}. Right: initial conditions \eqref{ini2}.}
\label{fig1}
\end{figure}

 In Fig. \ref{fig1} we have different numerical graphics of conformal trajectories of $V$. In all them, we have solved the system \eqref{eq2}, which in coordinates writes as 
\begin{equation}
\left\{ \begin{split}
  x''(s)&=q(y(s) z'(s)-y'(s) z(s))\\
  y''(s)&=q(x'(s) z(s)-x(s) z'(s))\\
  z''(s)&=q(x(s) y'(s)-x'(s) y(s)),
  \end{split}
  \right.
  \end{equation}
where  $\gamma(s)=(x(s),y(s),z(s))$. We have considered  different values of $q$ and initial conditions. In Fig. \ref{fig1}, left,  the initial conditions   are
\begin{equation}\label{ini1}
\gamma(0)=(1,0,0),\quad \gamma'(0)=(0,1,0).
\end{equation}
The second condition on $\gamma'(0)$ implies that $\gamma$ is parametrized by arc-length. Since $a_0=\langle \gamma(0),\gamma'(0)\rangle$ by \eqref{3}, we obtain $a_0=0$. By \eqref{4}, we know $a_1=1$. On the other hand, 
$|\gamma(0)\times\gamma'(0)|=1$. Thus $\kappa(s)=\kappa(0)=q|V(0)\times\gamma'(0)|=q$. 
Now $\kappa(s)=q$  and $\tau(s)=qs$.  In Fig \ref{fig1}, right, the initial conditions are 
\begin{equation}\label{ini2}
\gamma(0)=(1,0,0),\quad \gamma'(0)=(\frac{\sqrt{2}}{2},\frac{\sqrt{2}}{2},0).
\end{equation}
Then $a_0=\langle \gamma(0),\gamma'(0)\rangle=\frac{\sqrt{2}}{2}$. Thus $\kappa(s)=\sqrt{2} q/2$ and $\tau(s)=q(s+\frac{\sqrt{2}}{2})$. 

\section{The unit sphere $\s^3$}\label{s4}
 
In this section we study the conformal trajectories on the unit sphere $\s^3$ for all proper conformal vector fields.  We know from Sect. \ref{s2} that the proper conformal vector fields on $\s^3$ are generated by  the vector fields $V$ defined in \eqref{v2}. Recall that    $\lambda(p)=-2\langle \vec{a},p\rangle$.

\begin{theorem}\label{t2}
Let $V\in\mathfrak{X}(\s^3)$ be a proper conformal vector field given by \eqref{v2}. 
Let $\gamma\colon I\to\s^3$ be a curve parametrized by arc-length. If $\gamma$ is  a conformal trajectory of $V$, then 
\begin{equation}\label{ss} 
\langle \gamma(s),\vec{a}\rangle=a_1\cos(s)+a_2\sin(s),  
\end{equation}
 where $a_1,a_2\in\r$. If $\gamma$ is a geodesic, then $\gamma$ is a great circle of 
 $\s^3$ obtained by intersecting the sphere with a 2-plane containing $\vec{a}$. 
 Otherwise, the curvature of $\gamma$ is constant, 
 $\kappa(s)=|q|\sqrt{|\vec{a}|^2-(a_1^2+a_2^2)}$, 
 and the torsion is 
$$
\tau(s)= -q(a_1\sin(s)-a_2\cos(s)).
$$
 \end{theorem}
 
 \begin{proof}
 
 We know the relation between the Levi-Civita connection  $\nabla$ in $\s^3$ and 
 that of $\r^4$, $\nabla^e$, namely,
\begin{equation*}
\nabla_XY=\nabla_X^eY+\langle X,Y\rangle p.
\end{equation*}
As $\gamma$ is parametrized by arc-length, we deduce that 
$\nabla_{\gamma'}\gamma'=\gamma''+\gamma$. 

Let $\omega$ be the volume element on $\s^3$ defined by $\omega(u,v,w)=\mbox{det}(u,v,w,p)$, where $u,v,w\in T_p\s^3$. If $u,v\in T_p\s^3$,   
the cross product $u\times v$ on $\s^3$ is defined as the unique vector in $T_p\s^3$ such that $\langle u\times v,w\rangle=\omega (u,v,w)=\mbox{det}(u,v,w,p)$ for any $w\in T_p\s^3$. Therefore,  the cross product $V\times \gamma'$ may be expressed as
\begin{equation}
\label{cross-r4}
V\times\gamma'=\left|\begin{matrix}
a&x&x'&\partial_{x}\\
b&y&y'&\partial_{y}\\
c&z&z'&\partial_{z}\\
d&t&t'&\partial_{t}
\end{matrix}\right|
\end{equation}
where $\vec{a}=(a,b,c,d)$.
Obviously, $V\times \gamma'$ is perpendicular to $\vec{a}$. 
Thus, taking the scalar product with $\vec{a}$ in \eqref{eq1} we obtain 
$$
\frac{d^2}{ds^2}\langle \gamma,\vec{a}\rangle+\langle \gamma,\vec{a}\rangle=0.
$$
It follows that there exist two constants $a_1$ and $a_2$ such that
$$
\langle \gamma(s),\vec{a}\rangle=a_1\cos(s)+a_2\sin(s).
$$
By definition of curvature, we have  
$\kappa^2=\langle \nabla_{\gamma'}\gamma',\nabla_{\gamma'}\gamma'\rangle$. 

The special case when $\gamma$ is a geodesic   implies on one hand that 
$\gamma''(s)+\gamma(s)=0$ and on the other hand $V(s)\times \gamma'(s)=0$. 
The first equation yields $\gamma(s)=\cos(s) \vec{v}+\sin(s)\vec{w}$, 
where $\vec{v}$ and $\vec{w}$ are unitary and perpendicular. The second condition
implies that the first three columns in the determinant \eqref{cross-r4} are
linearly dependent. As a consequence, $\vec{a}$ lies in the 2-plane spanned by 
$\vec{v}$ and $\vec{w}$. This proves the first statement of theorem.

If $\gamma$ is not a geodesic, the unit normal vector is defined by
$\nabla_{\gamma'}\gamma'=\kappa N$. 
To obtain the curvature, we use \eqref{ss} to compute
\begin{equation*}
\begin{split}
\langle V\times\gamma',V\times\gamma'\rangle&=\langle V,V\rangle \langle \gamma',\gamma'\rangle-\langle V,\gamma'\rangle^2=\langle\vec{a},\vec{a}\rangle-\langle\gamma,\vec{a}\rangle^2-\langle\gamma',\vec{a}\rangle^2\\
&=|\vec{a}|^2-(a_1^2+a_2^2).
\end{split}
\end{equation*}
Hence $\gamma$ has constant curvature 
$\kappa=|q|\sqrt{|\vec{a}|^2-(a_1^2+a_2^2)}$. The expression of the torsion $\tau$ is obtained by Thm. \ref{t0} and identity \eqref{ss}. 
 
\end{proof}
 
 We show an explicit example of non-geodesic conformal trajectory.
 
 \begin{example}
 We prove that non-great circles of $\s^3$ are conformal trajectories. To be precise, let $\vec{a}=(a,b,c,d)\in\r^4$ and the corresponding conformal vector field $V(p)=\vec{a}-\langle\vec{a},p\rangle p$.   In the 3-dimensional space 
 $\vec{a}^\perp$ consider the orthonormal basis $\{\epsilon_1,\epsilon_2,\epsilon_3\}$ given by 
 $$
 \epsilon_1=\frac{(-d,c,-b,a)}{|\vec{a}|}~,
  \epsilon_2=\frac{(-c,-d,a,b)}{|\vec{a}|}~,
   \epsilon_3=\frac{(-b,a,d,-c)}{|\vec{a}|}~.
 $$
 Let $R\in(0,1)$ and consider the curve $\gamma$ on $\s^3$ defined by
 $$
 \gamma(s)=R\Big(\cos\frac{s}{R}~\epsilon_1+\sin\frac{s}{R}~\epsilon_2\Big)+
 \sqrt{1-R^2}~\epsilon_3.
 $$
 Obviously $\gamma$ is parametrized by the arc-length  and  $V(\gamma(s))=\vec{a}$ for any $s$.  Straightforward computations yield
 \begin{equation*}
\begin{split}
\nabla_{\gamma'}\gamma'&=\gamma''(s)+\gamma(s)=\frac{R^2-1}{R}\left(\cos\frac{s}{R}~\epsilon_1
 +\sin\frac{s}{R}~\epsilon_2\right)+ \sqrt{1-R^2}~\epsilon_3,\\
 V(\gamma(s))\times\gamma'(s)&=|\vec{a}|\left(\sqrt{1-R^2}\cos\frac{s}{R}~\epsilon_1
 +\sqrt{1-R^2}\sin\frac{s}{R}~\epsilon_2-R ~\epsilon_3\right).
 \end{split}
\end{equation*}
By choosing $q=-\frac{\sqrt{1-R^2}}{|\vec{a}|R}$, then $\gamma$ satisfies Eq. \eqref{eq1}. Notice that 
$\langle\gamma(s),\vec{a}\rangle=0$ and thus $a_1=a_2=0$ in \eqref{ss}. This shows that the torsion of $\gamma$ is $0$, as expected for a  circle in $\s^3$.

 \end{example}

\begin{figure}[hbtp] 
\begin{center}
\includegraphics[width=.25\textwidth]{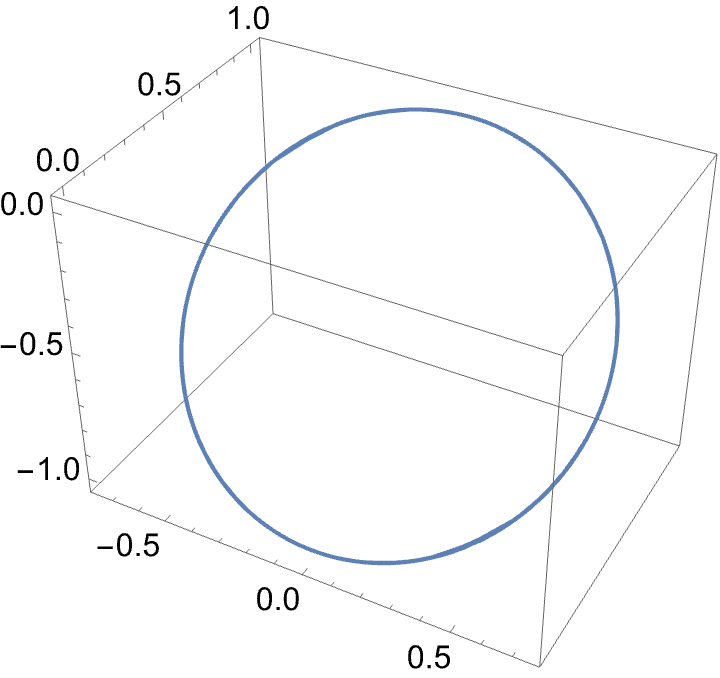}\quad\includegraphics[width=.32\textwidth]{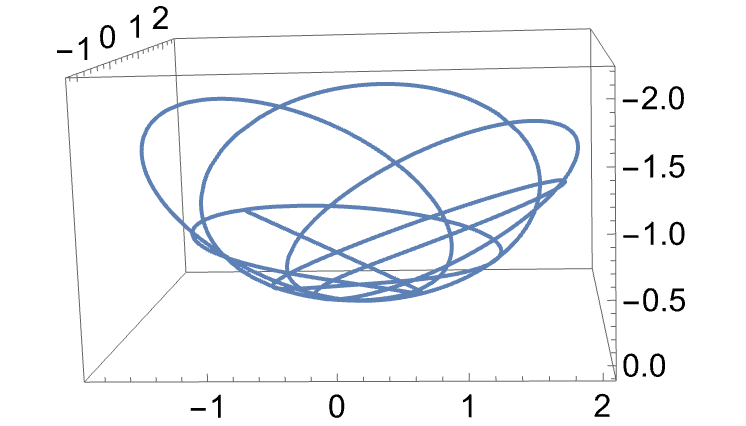}\quad\includegraphics[width=.32\textwidth]{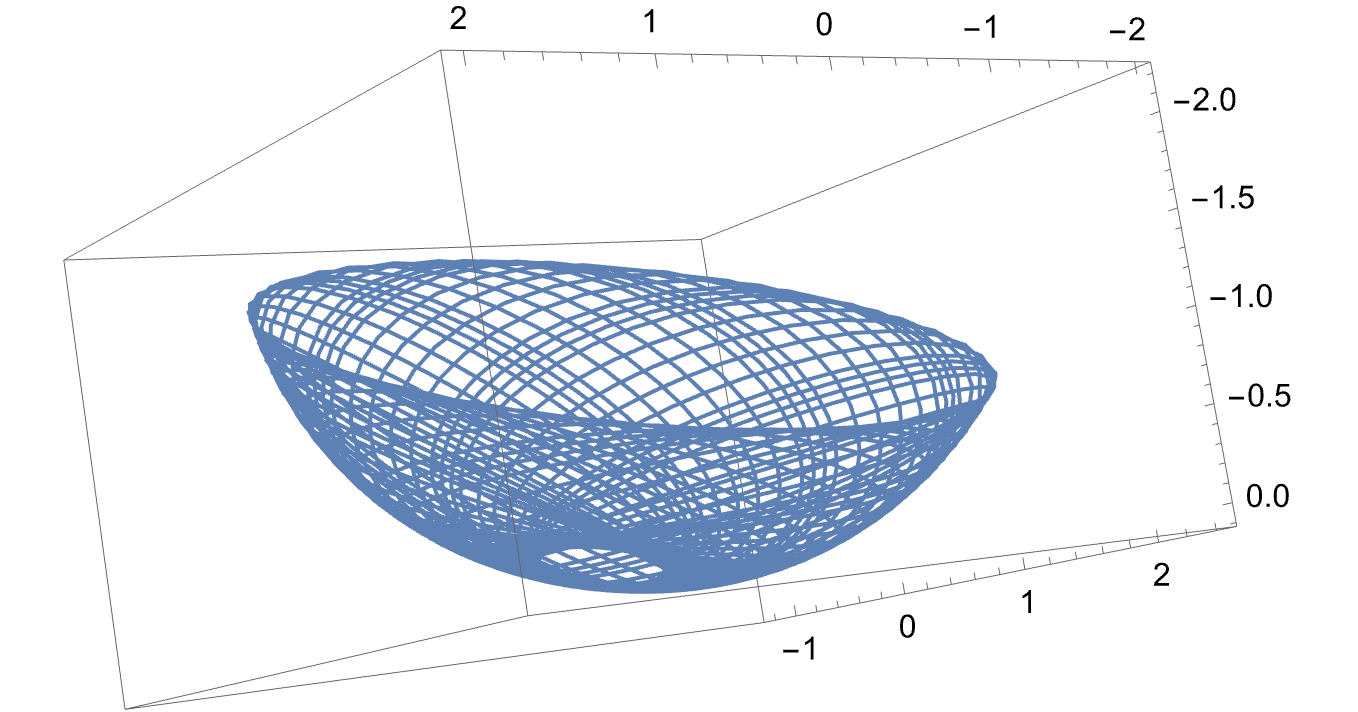}
\end{center}
\caption{Conformal trajectories of the vector field \eqref{v2} when $q=1$ and  $\gamma(0)=(0,1,0,0)$. Left: $\gamma'(0)=(-1,0,0,0)$. Middle and right: $\gamma'(0)=(0,0,1,1)$.}
\label{fig2}
\end{figure}

 We particularize for a  conformal vector field. Let   $\vec{a}=(0,0,0,1)$ in \eqref{v2}. Then 
 $$V=-tx\partial_x-t y\partial_y-t z\partial_z+(1-t^2)\partial_t.$$
 Here   $\lambda(x,y,z,t)=-2t$. If $\gamma(s)=(x(s),y(s),z(s),t(s))$ is parametrized by arc-length, then $\langle\gamma(s),\vec{a}\rangle=t(s)$ and $t(s)$ satisfies the differential equation  
$$t''(s)+t(s)=0.$$
In Fig. \ref{fig2} we show some graphics of conformal trajectories of $V$ once we take stereographic projection from $\s^3$ in $\r^3$.

 \section{The hyperbolic space $\h^3$} \label{s5}
 Consider the Lorentzian model  for $\h^3$. Here we follow the notation of  Sect. \ref{s2}. We know  that the proper conformal vector fields on $\h^3$ are generated by  the vector field  defined in \eqref{v3}.

\begin{theorem}\label{t3}
Let $V\in\mathfrak{X}(\h^3)$ be a proper conformal vector field given by \eqref{v3}. 
Let $\gamma\colon I\to\h^3$ be a curve parametrized by arc-length. If $\gamma$ is  a conformal trajectory of $V$, then 
\begin{equation}\label{ss2} 
\langle \gamma(s),\vec{a}\rangle=a_1\cosh(s)+a_2\sinh(s),  
\end{equation}
 where $a_1,a_2\in\r$. If $\gamma$ is a geodesic, then $\gamma$ is  obtained by intersecting $\h^3$ with a 
 Minkowski $2$-plane containing $\vec{a}$.  If $\gamma$ is not a geodesic, then  the curvature of $\gamma$ is constant, 
 $\kappa(s)=|q|\sqrt{| |\vec{a}|^2+a_1^2-a_2^2|}$, 
 and the torsion is 
$$
\tau(s)= q(a_1\sinh(s)+a_2\cosh(s)).
$$
 \end{theorem}
 
 \begin{proof}
 The proof is similar as in Thm. \ref{t2}.  We know the relation between the Levi-Civita connection  $\nabla$ in $\h^3$ and 
 that of $\r^4$, $\nabla^e$, namely,
\begin{equation*}
\nabla_XY=\nabla_X^eY-\langle X,Y\rangle p.
\end{equation*}
As $\gamma$ is parametrized by arc-length, we deduce that 
$\nabla_{\gamma'}\gamma'=\gamma''-\gamma$. Since $\langle V\times\gamma',\vec{a}\rangle=0$, from \eqref{eq1} we have  
$$
\frac{d^2}{ds^2}\langle \gamma,\vec{a}\rangle-\langle \gamma,\vec{a}\rangle=0.
$$
This gives \eqref{ss2}. 
If  $\gamma$ is a geodesic  on $\h^3$, then  $\gamma''(s)-\gamma(s)=0$ and this implies 
$\gamma(s)=\cosh(s) \vec{v}+\sinh(s)\vec{w}$, where $\langle\vec{v},\vec{v}\rangle=-1$, 
$\langle\vec{w},\vec{w}\rangle=1$ and $\langle\vec{v},\vec{w}\rangle=0$.  
Since    $V(s)\times \gamma'(s)=0$, then  $\vec{a}$ lies in the 2-plane spanned by 
$\vec{v}$ and $\vec{w}$.

Suppose now that  $\gamma$ is not a geodesic. Using \eqref{ss2} and $\kappa N=  V\times\gamma'$, we obtain
\begin{equation*}
\begin{split}
\kappa^2&=\langle V\times\gamma',V\times\gamma'\rangle=\langle V,V\rangle \langle \gamma',\gamma'\rangle-\langle V,\gamma'\rangle^2\\
&=\langle\vec{a},\vec{a}\rangle+\langle\gamma,\vec{a}\rangle^2-\langle\gamma',\vec{a}\rangle^2=|\vec{a}|^2+ a_1^2-a_2^2.
\end{split}
\end{equation*}
Thus $\gamma$ has constant curvature. For the last part of the theorem, we use  Thm. \ref{t0} in combination with \eqref{ss2}.

\end{proof}

\begin{example}
Let $\vec{a}=(1,0,0,0)$ and $\alpha>0$. Consider the curve $\gamma:\r\rightarrow\h^3$ given by
$$
\gamma(s)=
 \left(0,s,\alpha s^2+\alpha-\frac{1}{4\alpha},\alpha s^2+\alpha+\frac{1}{4\alpha}\right).
$$
Obviously $\gamma$ is parametrized by the arc-length and $V(\gamma(s))=\vec{a}$ for any $s$. Direct calculations lead to $\gamma''-\gamma=V\times \gamma'$ meaning that $\gamma$ is a 
(non-geodesic) conformal trajectory on $\h^3$ with $q=1$. Moreover, 
$\langle\gamma''-\gamma,\gamma''-\gamma\rangle=1$ which implies that the curvature of $\gamma$ is
$\kappa=1$. Then we obtain $\tau=0$ and the Frenet frame is
$$\left\{
\begin{array}{l}
T(s)=(0,1,2\alpha s,2\alpha s),\\[2mm]
N(s)=\left(0,-s,-\alpha s^2+\alpha+\frac{1}{4\alpha},-\alpha s^2+\alpha-\frac{1}{4\alpha}\right),
\\[2mm]
B(s)=(1,0,0,0).
\end{array}
\right.
$$
\end{example}
\begin{example}
Let $\vec{a}=(1,0,0,0)$ and $r>0$. Consider the curve 
$\gamma:(-\frac{\pi}{2},\frac{3\pi}{2})\rightarrow\h^3$ defined by
$$
\gamma(s)=\left(0,\frac{\cos(\phi(s))}{1+\sin(\phi(s))},
r-\frac{1}{2r(1+\sin(\phi[s)))},r+\frac{1}{2r(1+\sin(\phi[s)))}\right),
$$
where $\phi$ is a smooth function satisfying the differential equation $\phi'(s)=1+\sin(\phi(s))$.
It is straightforward to show that $\gamma$ is a unit speed conformal trajectory 
with $q=1$ having constant curvature $\kappa=1$ and zero torsion.
\end{example}

\begin{proposition}
If $\gamma$ is conformal trajectory on $\h^3$ with zero torsion, 
then $\vec{a}$ must be a space-like vector.
\end{proposition}
\begin{proof}
The condition $\tau=0$ yields $a_1=a_2=0$; see Thm. \ref{t3}. Moreover, it follows
that $\langle\vec{a},\gamma(s)\rangle=0$, for any $s$. Making use of the 
Cauchy–Bunyakovsky–Schwarz inequality (in $\mathbb{R}^3$), we obtain that 
$\langle\vec{a},\vec{a}\rangle\geq0$ 
and $\langle\vec{a},\vec{a}\rangle=0$ if and only if $\vec{a}=\vec{0}$.
\end{proof}
 
\begin{figure}[hbtp] 
\begin{center}
\includegraphics[width=.3\textwidth]{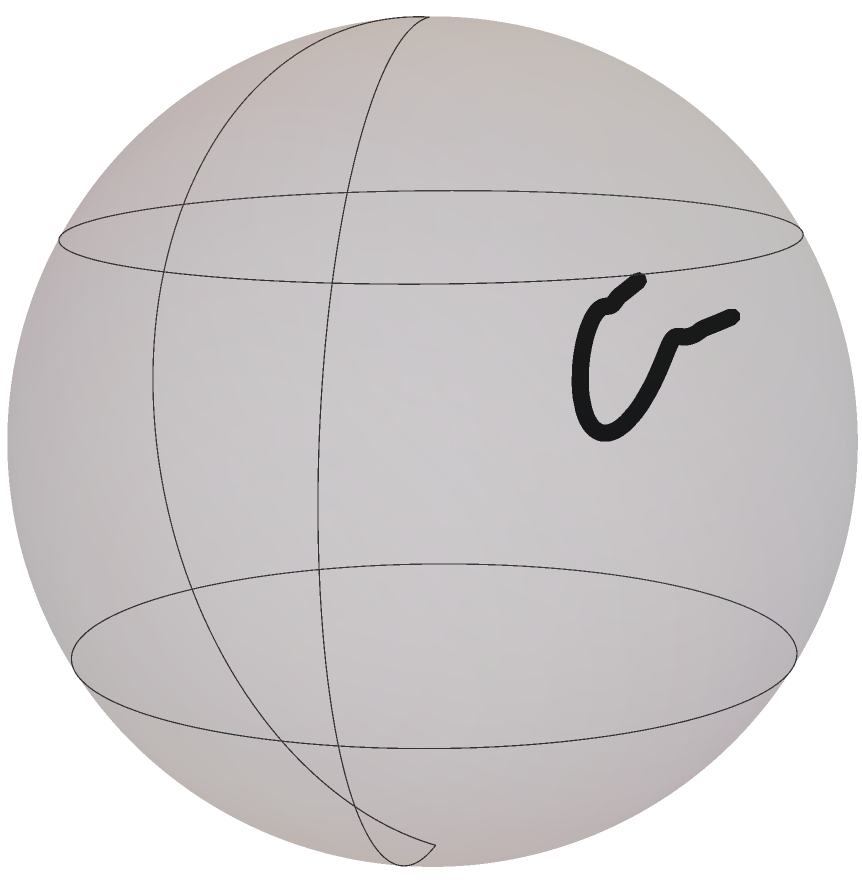}\qquad\includegraphics[width=.4\textwidth]{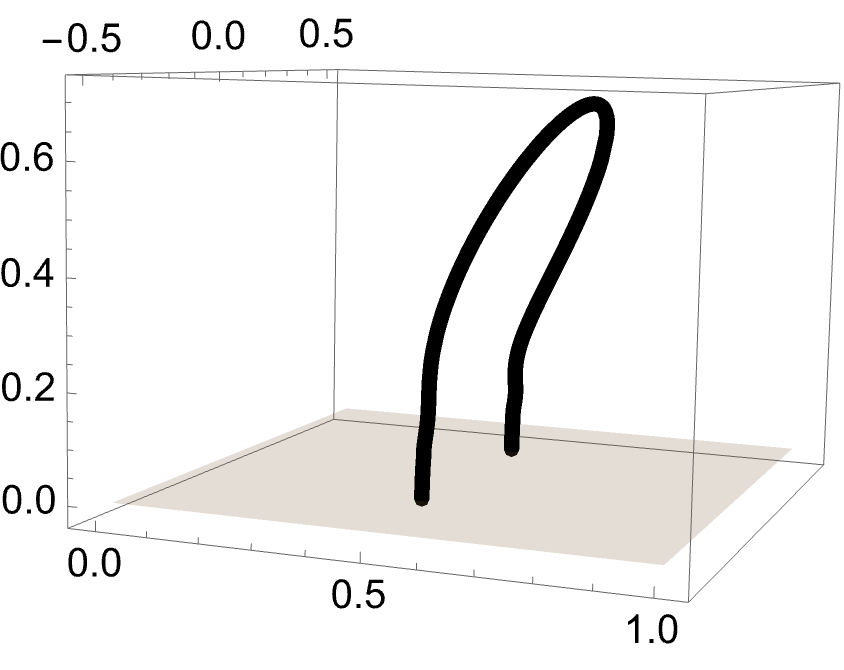}
\end{center}
\caption{A conformal trajectory for $V$  with $\vec{a}=(1,0,0,0)$ and initial conditions $\gamma(0)=(1,0,0,\sqrt{2})$, $\gamma'(0)$. Here $q=1$. Left: in the model of the Poincar\'e ball. Right: in the upper half-space model. }
\label{fig3}
\end{figure}

To show graphics of conformal trajectories, let us fix   $\vec{a}=(1,0,0,0)$ in \eqref{v3}. Then 
 $$V=(1+x^2)\partial_x+xy\partial_y+xz\partial_z+xt\partial_t.$$ 
 If $\gamma(s)=(x(s),y(s),z(s),t(s))$ and fix $q=1$, then  \eqref{eq1} writes as
 \begin{equation}\label{sys}
 \left.
 \begin{split}
 x''-x&=0\\
 y''-y-zt'+tz'&=0\\
 z''-z-ty'+yt'&=0\\
 t''-t+yz'-zy'&=0.
 \end{split}\right\}
 \end{equation}
In Fig. \ref{fig3} we show the  graphic of a conformal trajectory $\gamma$ of $V$ in the
 Poincar\'e ball (left) and in the upper half-space model (right).  
 In order to have a curve with non-zero torsion, we need $\langle\gamma(s),\vec{a}\rangle $ 
 is not $0$ at the initial condition. For this, it is enough $x(0)\not=0$. 
 Let $\gamma(0)=(1,0,0,\sqrt{2})$. Since in the system \eqref{sys} the curve $\gamma$ is
  parametrized by arc-length, then we take $\gamma'(0)$ of modulus $1$, for example 
  $\gamma'(0)=(0,1,0,0)$.  

\section*{Acknowledgements}  

 The authors  would like to thank Prof. Sebasti\'an Montiel (Granada)  for discussions on conformal vector fields in space forms.   Rafael L\'opez  is a member of the IMAG and of the Research Group ``Problemas variacionales en geometr\'{\i}a'',  Junta de Andaluc\'{\i}a (FQM 325). This research has been partially supported by MINECO/MICINN/FEDER grant no. PID2020-117868GB-I00, and by the ``Mar\'{\i}a de Maeztu'' Excellence Unit IMAG, reference CEX2020-001105- M, funded by MCINN/AEI/10.13039/501100011033/ CEX2020-001105-M.
Marian Ioan Munteanu is thankful to Romanian Ministry of Research, Innovation and Digitization, within Program 1 – Development of the national RD system, Subprogram 1.2 – Institutional Performance – RDI excellence funding projects, Contract no.11PFE/30.12.2021, for financial support.

\end{document}